\documentclass[12pt]{amsart}
\usepackage{amssymb}
\usepackage{amsmath}
\usepackage{pictex}
\newcommand\AND{\quad\mbox{and}\quad}
\newcommand\Area{\mbox{\sl Area}}
\newcommand\cog{\mbox{\sl cog}}
\newcommand\df{\mathfrak d}
\newcommand\ep{\varepsilon}
\newcommand\F{\mathbb F}
\newcommand\Ga{\Gamma}
\newcommand\id{\mbox{\sl id}}
\newcommand\la{\lambda}
\newcommand\lras{\buildrel * \over \leftrightarrow}
\newcommand\R{\mathbb R}
\newcommand\tos{\buildrel * \over \to}
\newcommand\Vol{\mbox{\sl Vol}}

\numberwithin{equation}{section}

\newtheoremstyle{mythm}
  {9pt}
  {9pt}
  {\itshape}
  {0pt}
  {\bfseries}
  {}
  { }
  {\thmnumber{(#2)}\thmname{ #1}\thmnote{ #3}}

\newtheoremstyle{mydef}
  {9pt}
  {9pt}
  {\normalfont}
  {0pt}
  {\bfseries}
  {}
  { }
  {\thmnumber{(#2)}\thmname{ #1}\thmnote{ #3}}

\theoremstyle{mythm}
\newtheorem{thm}[equation]{Theorem.}
\newtheorem{pro}[equation]{Proposition.}
\newtheorem{lem}[equation]{Lemma.}

\theoremstyle{mydef}
\newtheorem{dfn}[equation]{Definition.}
\newtheorem{exa}[equation]{Example.}

\begin{document}
\title{\large Non-backtracking random walks and cogrowth of graphs}
\author{\bf Ronald ORTNER and Wolfgang WOESS}
\email{woess@TUGraz.at,  ortner@finanz.math.tu-graz.ac.at}
\address{\parbox{1.4\linewidth}{Institut f\"ur Mathematik C, 
Technische Universit\"at Graz\\
Steyrergasse 30, A-8010 Graz, Austria}}
\date{June 27, 2003}
\thanks{Supported by FWF (Austrian Science Fund) project P15577}
\begin{abstract}
Let $X$ be a locally finite, connected graph without vertices of degree
$1$. Non-backtracking random walk moves at each step with equal
probability to one of the ``forward'' neighbours of the actual state,
i.e., it does not go back along the preceding edge to the preceding state. 
This is not a Markov chain, but can be turned into a Markov chain whose 
state space is the set of oriented edges of $X$. Thus we obtain for
infinite $X$ that the $n$-step non-backtracking transition probabilities tend to
zero, and we can also compute their limit when $X$ is finite.
This provides a short proof of old results concerning cogrowth of
groups, and makes the extension of that result to arbitrary regular
graphs rigorous.
Even when  $X$ is non-regular, but \emph{small cycles are dense in} $X$,
we show that the graph $X$ is non-amenable if and only if the non-backtracking 
$n$-step transition probabilities decay exponentially fast. 
This is a partial generalization
of the cogrowth criterion for regular graphs which comprises the original
cogrowth criterion for finitely generated groups of 
{\sc Grigorchuk} and {\sc Cohen}. 
\end{abstract} 
\maketitle

\markboth{\sf R. Ortner and W. Woess}{\sf Non-backtracking random walks and
cogrowth}
\baselineskip 15pt

\section{Introduction and results}\label{intro}
Let $X$ be the vertex set of a locally finite, connected graph, 
possibly with multiple edges and loops.
We write $e(x,y)$ for the number of edges between the vertices
$x$ and $y$,  if $y \ne x$, while $e(x,x)$ is \emph{twice} 
the number of loops at $x$ (see \S \ref{cogrowth}.B for a discussion).
The degree of a vertex $x \in X$ is $\deg(x)=\sum_y e(x,y)$.
We assume that $\deg(x) \ge 2$ for all $x \in X$.
\emph{Non-backtracking (simple) random walk} (NBRW) is the following 
random process: at the beginning, the walker starts at some vertex $x$
and chooses with equal probablity one of the incident edges. He steps
to the other end of that edge. At the later steps, the rule is the same,
but the walker selects with equal probability only among those incident
edges that are different from the one transversed at the previous step.

We write $q^{(n)}(x,y)$ for the probability that the random walker,
starting at vertex $x$, is at vertex $y$ at the $n$-th step.
Note that NBRW is \emph{not} a Markov chain on $X$. The defining
property of a Markov chain, that ``the future depends only on the
actual state and not on the past'', is violated, since the walker has to
remember the edge along which he reached the actual state before 
moving on. 

However, it is easy to turn NBRW into a Markov chain by
changing the state space: with each edge, we associate two oppositely
oriented edges $e, \check e$ (with $\check{\check e} = e$). We write $e^-$ and $e^+$
for the initial and terminal vertex of the edge $e$, so that
$(\check e)^- = e^+$ and $(\check e)^+ = e^-$. (Note in particular, that
for each a priori unoriented loop we get two oriented ones !)
We now consider NBRW as a \emph{Markov} process whose new state 
space is the set $E = E(X)$
of oriented edges, with transition matrix 
$Q_E = \bigl(q_E(e,f)\bigr)_{e,f\in E}$ given by
$$
q_E(e,f) = \begin{cases} \dfrac{1}{\deg(e^+)-1}\,,
                    &\mbox{if $e \to f$, that is, $f^- = e^+$ and $f \ne \check e$\,,}\\
		 0\,,&\mbox{otherwise.}   
         \end{cases}
$$
Then there is the following obvious link between edge-NBRW and 
vertex-NBRW: for vertices $x, y \in X$,
\begin{equation}\label{qn}
q^{(n)}(x,y) = \frac{1}{\deg(x)} 
\sum_{\scriptstyle e,f \in E : \atop \scriptstyle e^+=x,f^+=y}
q_E^{(n)}(e,f)\,,
\end{equation}
where $q_E^{(n)}$ denotes the $n$-step transition probabilities, i.e.,
the elements of the matrix power $Q_E^n$, with $Q_E^0 = I_E$, the identity
matrix over $E$. (Attention: $q^{(n)}(x,y)$ is \emph{not} the $(x,y)$-element
of an $n$-th matrix power over $X$\,!)

The following result is then a consequence of basic
Markov chain theory.

\begin{thm}\label{limit1}
{\rm (a)} If $X$ is finite, connected, with minimum degree $2$, 
then for all $x,y \in X$,  
$$
\lim_{n \to \infty} \frac{1}{n}\Bigl( q^{(1)}(x,y) + q^{(2)}(x,y)
+ \dots + q^{(n)}(x,y)\Bigr) = \frac{\deg(y)}
{|E(X)|}\,.
$$

\smallskip\noindent
{\rm (b)} If in addition to the assumptions of \/{\rm (a),} $X$ has 
minimum degree $3$, then for all $x,y \in X$,
$$
\begin{aligned}
\lim_{n \to \infty} q^{(n)}(x,y) &= \frac{\deg(y)}{|E(X)|}\,,\quad
\mbox{if $X$ is not bipartite,}\\
\lim_{n \to \infty} q^{(2n+\delta)}(x,y) &= \frac{2\deg(y)}{|E(X)|}\,,\quad
\mbox{if $X$ is bipartite,}
\end{aligned}
$$
where $\delta=0$ (resp. $\delta=1$) according to whether $x$ and $y$ are 
at even (resp. odd) distance.

\medskip\noindent
{\rm (c)} If $X$ is infinite, connected, with minimum degree $2$, 
then for all $x, y \in X$, 
$$
\lim_{n \to \infty} q^{(n)}(x,y) = 0\,.
$$
\end{thm}

In statements (a) and (b), note that $|E(X)|$ is twice the number of 
non-oriented edges.

As usual, the \emph{distance} $d(x,y)$ between two vertices $x,y \in X$
is the minimum length of a path connecting the two. The
\emph{ball} of radius $R$ centred at $x$ is the subgraph
$B(x,R) = \{ y \in X : d(y,x) \le R \}$ of $X$. 
Recall that a \emph{cycle} of length $n$ in $X$ consists of a sequence
$e_n = e_0, \dots, e_{n-1}$ of distinct edges whose initial vertices
are all distinct, such that $e_{k-1} \to e_k$
for all $k=1, \dots, n$.

\begin{dfn}\label{cycles} We say that \emph{small cycles
are dense in}  $X$, if there is $R > 0$ such that every ball $B(x,R)$
in $X$ contains a cycle.
\end{dfn}

Every finite, connected graph with minimum degree $2$ 
satisfies this condition.

The \emph{automorphism group} of $X$ consists of all bijections 
$g: X \to X$ which satisfy $e(gx,gy) = e(x,y)$ for all $x,y \in X$.
A graph is called \emph{transitive}, resp.\ \emph{almost transitive} 
if the automorphism group acts with one orbit, resp. finitely
many orbits on $X$. Obviously, an infinite, almost transitive graph
with minimum degree $2$ has dense small cycles unless it is a tree.
(To be precise, we require of a tree that it does not have
multiple edges.)

\begin{lem}\label{radius}
If small cycles are dense in $X$ then 
$$
\rho(Q) = \limsup_{n \to \infty} q^{(n)}(x,y)^{1/n}
$$
is independent of $x, y \in X$, and $0 < \rho(Q) \le 1$.
(If $X$ is finite then $\rho(Q) = 1$.)
\end{lem} 

The following strengthens Theorem \ref{limit1}\,(c) for
almost transitive graphs.

\begin{thm}\label{limit2}
If $X$ is infinite, connected, with minimum degree $2$, and
almost transitive, then for all $x, y \in X$, 
$$
\lim_{n \to \infty} q^{(n)}(x,y)/\rho(Q)^n = 0\,.
$$
\end{thm}

The \emph{isoperimetric constant} $\iota(X)$ of a connected, locally 
finite graph $X$ is 
$$
\iota(X) = \inf \left\{ \frac{\Area(F)}{\Vol(F)} : 
                        F \subset X \;\mbox{finite} \right\}\,,
$$
where $\Vol(F) = \sum_{x \in F} \deg(x)$ and $\Area(F)$ is the
number of edges with one endpoint in $F$ and the other in 
$X \setminus F$. The graph is called \emph{amenable} if $\iota(X) = 0$.
Non-amenable graphs are also called (infinite) \emph{expanders}.

Consider the Hilbert space $\ell^2(E)$ of all functions $F: E \to \R$ with 
$\langle F,F \rangle < \infty$, with the ordinary inner product
$$
\langle F,G \rangle = \sum_{e \in E} F(e)G(e).
$$
Then $Q_E$ acts on this space by $Q_EF(e) = \sum_{f \in E} q_E(e,f)F(f)$.
We denote by $\|Q_E\|$ the corresponding operator norm, and by 
$\rho_2(Q_E) = \lim_n \|Q_E^n\|^{1/n}$ its spectral radius. 
Note that $\rho(Q) \le \rho_2(Q_E) \le \|Q_E\|$ in general.

\begin{pro}\label{norm} {\rm (a)} One has always $\|Q_E\|=1$.\\[4pt]
{\rm (b)} If small cycles are dense in $X$, then $\rho(Q) = \rho_2(Q_E)$.
\end{pro}


\begin{thm}\label{amenable} Suppose that $X$ is connected, that small
cycles are dense, and that there is $M < \infty$ such that
$2 \le \deg(x) \le M$ for all $x \in X$.

\smallskip

Then $X$ is amenable if and only if $\rho(Q) = 1$.
%
\end{thm}

With these results and their proofs we aim principally at
extending and explaining previous material regarding
\emph{cogrowth} of graphs and groups and at shedding
new light on cogrowth by studying it in terms of NBRW on
the oriented edges. We also think that NBRW on the (oriented)
edge set of an arbitrary graph is an interesting random 
process in its own right.

In \S \ref{cogrowth}, we first recall (ordinary) simple random
walk on a graph and some of its basic properties in order to
put our results on NBRW in the right perspective.
We then consider cogrowth of  graphs, which is best understood
in terms of universal covering trees, and explain how
Theorems \ref{limit1}, \ref{limit2} and \ref{amenable} apply.
In \S 2 we also give various references.

\S \ref{proofs} is dedicated to the proofs of the results stated
here.

Some additional remarks and observations can be found in \S \ref{final}.

\section{Simple random walk, and cogrowth of graphs}\label{cogrowth}

\noindent
{\bf A. Simple random walk} (SRW) is mostly considered on graphs without
multiple edges, and loops are usually counted only once for the degree
of a vertex. Here, multiple edges are admitted, and
we count each loop twice. SRW is the Markov
chain on the (vertex set of the) graph $X$ with transition matrix
$P = \bigl( p(x,y) \bigr)_{x,y \in X}$ given by
$$
p(x,y) = \frac{e(x,y)}{\deg(x)}.
$$
Thus, contrary to NBRW, the walker does not remember from where he
did come at the previous step, and chooses at random any one among
the outgoing edges at the actual vertex. A possible interpretation
for counting each loop twice is that topologically, the walker
standing at a vertex $x$ sees two ``ends'' of each loop at $x$ among
which he may choose.
We write $p^{(n)}(x,y)$ for the $n$-step transition probability from $x$
to $y$.

The transition matrix $P$ acts by $Pg(x) = \sum_y p(x,y)g(y)$ on the
Hilbert space $\ell^2(X,\deg)$ of all functions $g: X \to \R$ with
$\langle g,g \rangle < \infty$, where the inner product is
$$
\langle g,h \rangle = \sum_{x \in X} g(x)h(x)\,\deg(x).
$$
We denote by $\|P\|$ the norm of this operator.

Here is a list of well-known properties of SRW. (Recall once more 
that $E = E(X)$ is the set of oriented edges as in \S \ref{intro}, so that
$|E(X)|$ is twice the number of ``ordinary'' non-oriented edges.)

\begin{pro}\label{SRW} Let $X$ be a connected, locally finite graph.
\begin{itemize}
\item[(a)] If $X$ is finite and not bipartite, then for all $x,y \in X$,
$$
\lim_{n \to \infty}  p^{(n)}(x,y) = \frac{\deg(x)}{|E(X)|}\,.
$$
If $X$ is finite and bipartite, then for all $x,y \in X$, with
$\delta \in \{0,1\}$ such that $d(x,y) \equiv \delta \mod 2$,
$$
\lim_{n \to \infty}  p^{(2n+\delta)}(x,y) = 2\frac{\deg(x)}{|E(X)|}\,.
$$
\item[(b)] If $X$ is infinite, then for all $x,y \in X$,
$$
\lim_{n \to \infty}  p^{(n)}(x,y) = 0\,.
$$
\item[(c)] The \emph{spectral radius}
$$
\rho(P) = \limsup_{n \to \infty} p^{(n)}(x,y)^{1/n}
$$
is independent of $x,y \in X$, and $\|P\| = \rho(P)$.
\smallskip
\item[(d)] If $X$ is infinite and almost transitive then
$$
\lim_{n \to \infty}  p^{(n)}(x,y)/\rho(P)^n = 0\,.
$$
\item[(e)] $X$ is amenable if and only if $\rho(P) =1\,.$
\end{itemize}
\end{pro}

Statements (a) and (b) follow from basic Markov chain theory,
see e.g.\ {\sc Chung} \cite{Chu} or {\sc Seneta} \cite{Sen}: 
the Markov chain given by $P$ is
\emph{irreducible} ($\forall\ x, y \in X \;\exists\;n=n(x,y)\ge 0$
such that $p^{(n)}(x,y) > 0$). Its \emph{period}
$\df(P) = \gcd \{ n : p^{(n)}(x,x) > 0 \}$ is $=2$ when $X$ is bipartite,
and $=1$, otherwise. Finally, $\mu(x) = \deg(x)$ defines an \emph{invariant}
measure. If $X$ is finite then $\mu(X)=|E(X)|$, and $\mu_0(x)=\mu(x)/|E(X)|$
is an invariant probability measure. Therefore, (a) follows from the basic
convergence theorem, see \cite{Chu}, Thm.\ 1 in \S I.6 or \cite{Sen}, Thm.\ 4.2. 
If $X$ is infinite then $\mu(X) = \infty$, whence
the random walk cannot be \emph{positive recurrent}, and (b) must hold.
We shall encounter these notions in more detail in \S \ref{proofs}.

For statement (c), see e.g.\ {\sc Woess} \cite{Wbook},  \S 10.
In particular, the fact that $\rho(P) = \rho_2(P)$, the $\ell^2$-spectral
radius of $P$, follows from self-adjointness of $P$ on $\ell^2(X,\deg)$.

Regarding statement (d), this is immediate when
$\sum_n  p^{(n)}(x,y)/\rho(P)^n < \infty\,$. If the series diverges
then it follows from  Theorem 7.8 in \cite{Wbook}
(which is basically due to {\sc Guivarc'h} \cite{Gui})
that $\rho(P)=1$, and we can apply (b).

Statement (e) has a long history, going back to {\sc Kesten}'s amenability
criterion for finitely generated groups \cite{Kes}. The version
stated here is due to {\sc Dodziuk and Kendall} \cite{Dod-Ken} based 
on a previous paper by {\sc Dodziuk} \cite{Dod}.


\medskip\noindent
{\bf B. Cogrowth} is a notion of asymptotic density of a graph.
It is best understood in terms of the \emph{universal cover} of
the graph $X$. This is a (unique) \emph{tree} $T$ together with a surjective
mapping $\pi: T \to X$ which is a local homeomorphism, i.e.,
if $\tilde x, \tilde y$ are neighbours in $T$ then so are 
$\pi(\tilde x), \pi(\tilde y)$ in $X$, and $\deg_T(\tilde x) =
\deg_X\bigl(\pi(\tilde x)\bigr)$ for every vertex $\tilde x \in T$.

The covering tree can be constructed as follows: a non-backtracking walk
of length $n\ge 0$ in $X$ is a sequence $e_1, \dots, e_n$ of edges such that 
$e_{k-1} \to e_{k}$
for $k=2, \dots, n$.
Its initial and terminal vertices are $e_1^-$ and $e_n^+$, respectively. 
If $n=0$, we have
an empty path, for which we have to specify its initial = terminal vertex.
We now choose a root (reference vertex) $o \in X$, and define $T$ as 
the set of all non-backtracking paths $\tilde x$ starting at $o$, 
including the empty path. Two such paths are defined to be neighbours
in $T$ if one of them extends the other by a single edge.
The mapping $\pi$ assigns to each $\tilde x \in T$ its terminal vertex
$x \in X$. 

Now let $x, y \in X$, and choose $\tilde x \in T$ such that $\pi(\tilde x)=x$.
Write $T(y) = \{ \tilde y \in T : \pi(\tilde y) = y\}$, and consider the sphere 
$S(\tilde x, n) = \{ \tilde v \in T : d_T(\tilde v,\tilde x) = n\}$,
where $d_T(\cdot,\cdot)$ is the distance in $T$.
Then \emph{(ordinary) cogrowth} at $x, y \in X$ is the sequence
\begin{equation}\label{ordinary}
\cog_n(x,y) = \frac{|S(\tilde x, n) \cap T(y)|}{|S(\tilde x, n)|}\,,\quad
n \ge 0\,.
\end{equation}
The graph $X$ being ``small'' corresponds to $\bigl( \cog_n(x,y) \bigr)_n$
being ``large''. Besides finiteness, also amenability is a ``smallness''
condition, whence it is natural to look for a link between cogrowth
and amenability.

Cogrowth was initially introduced by {\sc Grigorchuk} \cite{Gri} and
later {\sc Cohen} \cite{Coh} for finitely generated \emph{groups}. If
$\Ga$ is such a group, then we can represent it as a
factor $\F_s/N$, where $\F_s$ is the \emph{free group} on $s$ free
generators $\tilde a_1, \dots, \tilde a_s$, and $N$ is a normal subgroup
of $\F_s$. 
Let $\pi:\F_s \to \Ga$ be the factor map.
We write $\tilde a_{-i} = \tilde a_i^{-1}$ and set 
$\tilde S = \{ \tilde a_i : i = \pm 1, \dots, \pm s \}$.
Then the Cayley graph of $\F_s$ with respect to $\tilde S$ is the
$2s$-regular tree, which is the covering tree of the Cayley graph
of $\Ga$ with respect to the generators 
$a_i=\pi(\tilde a_i)$.
It is best to consider immediately the oriented edges of
that Cayley graph: every $x \in \Ga$ is the inital point of an edge
of type $\tilde a_i$, whose endpoint is $xa_i$; the associated ``inverse''
edge goes from $xa_i$ to $x$ and has type $\tilde a_{-i}$ 
($i = \pm 1, \dots, \pm s$). Every pair of this type corresponds to one
unoriented edge. Note that generators with $a_i=a_{-i} \ne \id$ give rise to 
multiple edges, and when $a_i=a_{-i} = \id$, we get loops. This also explains
why loops should be counted twice for the degrees. Thus, the factor map $\pi$
becomes the covering map from the tree onto the Cayley graph.

Note that for groups, $\cog_n(x,x)$ is the same for all $x$.
Amenability of a finitely generated group $\Ga$ is equivalent with
amenability of any of its (locally finite) Cayley graphs.
The main result of \cite{Gri} and \cite{Coh}, restated in our 
notation, was that 
\begin{equation}\label{classic}
\Ga \;\mbox{is amenable} \iff  \limsup_{n\to\infty} \cog_n(x,x)^{1/n} = 1\,.
\end{equation}

This has been generalized to regular graphs by {\sc Northshield} \cite{No1},
who was also the first to explain cogrowth in terms of covering trees.
One of the basic tools for studying cogrowth of regular graphs is a
functional equation between the generating functions
$C(x,y|t) = \sum_n \cog_n(x,y)\,t^n$ of the cogrowth sequence and 
$G(x,y|z) = \sum_n p^{(n)}(x,y)\,z^n$ of the transition probabilites of SRW:
if $X$ is $d$-regular then with our notation and normalizations,
\begin{equation}\label{functional}
C(x,y|t) =  \frac{1}{d}\delta_x(y) + 
\frac{(d-1)^2-t^2}{d(d-1+t^2)}G\bigl(x,y|z(t)\bigr)\,,\quad\mbox{where}\quad
z(t) =  \frac{dt}{d-1+t^2}\,,
\end{equation}
A first version of (\ref{functional}) is contained in the Ph.D. thesis of
{\sc Grigorchuk}.
Various proofs of that formula have appeared:
{\sc Woess} \cite{Wcog}, {\sc Szwarc} \cite{Szw} (both for groups), 
{\sc Northshield} \cite{No1} (shortest), {\sc Bartholdi} 
\cite{Bar} (more general). In spite of \cite{Bar}, there is no satisfactory
version of that formula for non-regular graphs. Nevertheless, {\sc Northshield} 
\cite{No2} proves a clever extension of (\ref{classic}) to \emph{quasi-regular}
graphs (non-regular graphs satisfying a certain uniform growth condition).
 
\smallskip

More generally, we can consider a sequence 
$\nu = (\nu_{\tilde x,n})_{\tilde x \in T, n \ge 0}$, where each 
$\nu_{\tilde x,n}$ is a probability measure concentrated on the sphere
$S(\tilde x,n)$ of radius $n$ centred at $\tilde x$ in the covering tree 
$T$ of $X$, with $\pi(\tilde x) = x$. 
Note that there is a natural bijection between $S(\tilde x, n)$ and
$S(\tilde x',n)$, when $\pi(\tilde x) = \pi(\tilde x')$. We require
that in this case, $\nu_{\tilde x',n}$ is the image of $\nu_{\tilde x,n}$
under that bijection. Then we can define
\begin{equation}\label{sequence}
\cog^{\,\nu}_n(x,y) = \nu_{\tilde x,n}\bigl(T(y)\bigr)\,,\quad x, y \in X\,,
\;\pi(\tilde x) = x\,.
\end{equation}
When each $\nu_{\tilde x,n}$ is equidistribution on $S(\tilde x,n)$,
this is ordinary cogrowth.

Another choice is to define
$$
\nu_{\tilde x,n}(\tilde y) = \frac{1}{\deg(\tilde x)}  
\frac{1}{\deg(\tilde x_1)-1} \cdots \frac{1}{\deg(\tilde x_{n-1})-1}\,,
$$
where $\tilde x, \tilde x_1, \dots, \tilde x_{n-1},\tilde y$ are the 
consecutive vertices on the unique path in $T$ from $\tilde x$ to 
$\tilde y \in S(\tilde x,n)$. Cogrowth with respect to this choice of $\nu$
is the same as NBRW:
\begin{equation}\label{NBcog}
\cog^{\,\nu}_n(x,y) = q^{(n)}(x,y)
\end{equation}
In the specific case of regular graphs, the two concepts coincide.
Thus, besides ordinary cogrowth, non-backtracking random walk is 
another way to extend cogrowth from regular to arbitrary graphs.

\section{Proofs}\label{proofs}

In this section, we always use the basic assumption that
$X$ is a locally finite, connected graph with minimum degree $2$.
 
It may be best to think of edge-NBRW as simple random walk on
the \emph{oriented line graph} (OLG) of $X$. This is the digraph
whose vertex set is $E=E(X)$, and there is an oriented (2nd order)
edge from $e$ to $f$ ($e,f \in E$) if $e \to f$.
Our Markov chain with transition matrix $Q_E$
is not symmetric, nor reversible like SRW on an unoriented graph.
However, the counting measure $\la$, given by
$\la(e) = 1$, is an invariant measure for $Q_E$, that is,
\begin{equation}\label{invariant}
\sum_{e \in E} \la(e)q_E(e,f) = \la(f) \quad \forall\ f \in E\,.
\end{equation}
We now recall a few basic Markov chain notions. We write $e \tos f$
if there is $n \ge 0$ such that $q_E^{(n)}(x,y) > 0$ (i.e., there is
an oriented path from $e$ to $f$ in the OLG, a transitive relation), 
and $e \lras f$ if
$e \tos f $ and $f \tos e$. The equivalence classes with respect to
the relation $\lras$ are called \emph{irreducible classes}. An 
\emph{essential} class $V$ is an irreducible class with the property 
that $e \in V$ and $e \tos f$ implies $f \in V$. Its elements are
also called essential. The Markov chain and its transition matrix $Q_E$
are called irreducible if the state space $E$ forms a single irreducible 
class. (In graph theoretic terminology, this means that the OLG is
\emph{strongly connected}.) 

\begin{lem}\label{irred} If  $X$ is finite then
$Q_E$ is irreducible, unless $X$ is a cycle.
\end{lem}

\begin{proof} Assume that $X$ is not a cycle.
Since $X$ is connected, for any pair of edges $e, f$, at least one of
$e \tos f$, $e \tos \check f$, $\check e \tos f$,
or $\check e \tos \check f$ must hold. Therefore it is sufficient to show
that $e \tos \check e$ for every $e \in E$.

Let us first assume that $e$ is not contained in any cycle of $X$.
As $\deg(x) \ge 2 \;\forall x$
we can find inductively a sequence $e=e_0, e_1, e_2, \dots$ of edges
such that $e_{k-1} \to e_k$.
By finiteness of $X$, there must be a minimal index $m$
such that $e_m^+ = e_i^-$ for some $i \in \{ 1, \dots, m-1 \}$.
The edges $e_i, \dots, e_m$ form a cycle $C_1$,
so that
$$
e =e_0 \tos e_m \to \check e_{i-1} \tos \check e_0 = \check e\,.
$$

Now assume that $e$ is contained in a cycle $C_1$ formed by edges
$e=e_0, \dots, e_m$. Since we are assuming that $X$ is not a cycle, there
is a vertex $e_i^-=:x$ in $C_1$ with $\deg(x) \ge 3$. Thus, there an edge
$f$ with $f^-=x$ such that $f\notin\{\check e_{i-1},e_i\}$ (for $i=0$ we
intend $e_{-1}=e_m$). If $f$ does not lie on any cycle in $X$, we have 
already seen that $f \tos \check f$, whence
$$
e =e_0 \tos e_{i-1} \to f \tos \check f \to  \check e_{i-1} \tos \check e_0 = \check e\,.
$$
On the other hand, assume that $f$ is contained in a cycle $C_2$ formed by edges
$f=f_0, \dots, f_\ell$. Then there must be another edge $f_k$ ($k>0$) incident with
some vertex in $C_1$. Let $j$ be the minimal index $\in \{ 0, \dots, m \}$
with $e_j^+=f_k^+$ for some $k \in \{ 1, \dots, \ell \}$. Then
\[
e \tos e_{i-1} \to f=f_0 \tos f_k \to  \check e_j \tos \check e\,.     \qedhere
\]
\end{proof}

If $X$ is a finite cycle, then the OLG consists of two disjoint, oriented
cycles of the same length, each of which constitutes an essential class
of $Q_E$, on which $NBRW$ moves ``forward'' deterministically.

\begin{lem}\label{infinite} If  $X$ is infinite then for any edge $e \in E$ there are
infinitely many edges $f\in E$ with $e \tos f$.
\end{lem}

\begin{proof}
Let $e \in E$ and $X'$ be the graph that results from $X$ by removing $e$
and $\check e$.
If $X'$ is connected then by infiniteness, $e \tos f$ for infinitely many 
$f \in E$. The same holds if $e$ is directed towards an infinite component. 
Thus, let us assume that $e$ is directed towards a finite component $X_1'$ 
of $X'$. By infiniteness of $X$, $\check e$ is directed to the other, infinite
component, so that $\check e \tos f$ for infinitely many $f \in E$. Applying
the method of proof of Lemma \ref{irred} to $X_1'$, we have $g \tos \check g$ 
for some edge $g$
with $e \to g$ in $X$ (remember that we assumed that $\deg(e^+)\geq 2$). 
It follows that $e \to g \tos \check g \to \check e$ and hence $e \tos f$ 
for infinitely many $f \in E$.
\end{proof}

In general, if $Q_E$ is irreducible, then we can define
its \emph{period} by
$$
\df = \df(Q_E) = \gcd \{ n : q_E^{(n)}(e,e) > 0 \}\,,
$$
which is independent of $e \in E$.

\begin{lem}\label{period}
Let $X$ be a finite, connected graph with $\deg(x)\ge 3$ for all $x\in X$. Then
the period of the associated edge-NBRW is either 2 or 1, depending on whether
$X$ is bipartite or not (respectively).
\end{lem}

\begin{proof}
First we shall show that $\df(Q_E)\in\{1,2\}$. Let $e,f,g$ be
three edges with $e^-=f^-=g^-=:x$. By Lemma \ref{irred} we have
$e \tos \check f$ and $e \tos \check g$. That is, there are
two non-backtracking closed paths at $x$ in $X$ both starting with $e$, one
terminating with $\check f$, the other one with $\check g$. Since the 
starting edge in both these paths is not the reversed terminating
edge, they can easily be turned into two cycles $C_1,C_2$ at $x$ 
formed by edges $e_1,\ldots,e_n$ and $f_1,\ldots,f_m$, respectively.
Both $C_1$ and $C_2$ start with the same edge $e_1=f_1=e$. We claim that 
we may assume that the second edges $e_2$ and $f_2$ in $C_1$ and $C_2$ (resp.)
do not coincide. Consider the case where $e_2=f_2$.
By assumption, $\deg(e^+)\geq 3$ and there is another edge 
$g\neq e_2,\check e$ with $e^+=g^-$. Since initial 
vertices do not occur more than once in each cycle, neither $C_1$ 
nor $C_2$ contains $g$. As $\deg(\cdot) \ge 3$ we can find 
inductively a sequence $g=g_1, g_2, g_3, \dots$ of edges with distinct 
initial vertices such that $g_{i} \to g_{i+1}$. By finiteness of $X$, 
there must be a minimal index $k$ such that $g_k^+$ occurs as initial 
vertex of an edge in one of the cycles $C_1,C_2$. Let us assume that 
$g_k^+=e_\ell^-$ for some $\ell \in \{ 1, \dots, n \}$. Then we 
may replace the cycle $C_1$ by $e_1,g_1,\ldots,g_k,e_\ell,\ldots,e_n$
so that the two cycles in $X$ have the claimed property. 
A similar argument shows
that we also may assume that $e_n\neq f_m$.

Thus we have two cycles of length $n$ and $m$, respectively.
Since we assumed $e_2\neq f_2$ and $e_n\neq f_m$, we have
$\check e_n \tos \check e_n$ in $n+m-2$ steps via
$$
\check e_n \to \ldots\to \check e_2 \to f_2 \to \ldots \to f_m \to \check e_n.
$$
Therefore, $\df(Q_E)$ must be a factor of $n$, $m$ and $n+m-2$,
whence $\df(Q_E) \in \{1,2\}$.


It is now clear that we must have $\df(Q_E) = 2$, if $X$ is bipartite.
Otherwise, $X$ contains an odd cycle, so that $q_E^{(k)}(e,e) >0$
for some odd $k$. Thus, we cannot have $\df(Q_E) = 2$, that is,
$\df(Q_E) = 1$.
\end{proof}

\begin{proof}[\bf Proof of Theorem \ref{limit1}]
(a+b) If $X$ is finite, but not a cycle, then we can use Lemma \ref{irred}.
Let $e, f \in E$ and $r \ge 0$ such that $q_E^{(r)}(e,f) > 0$.
Then $q_E^{(n)}(e,f) > 0$ if and only if $n \equiv r \mod \df$ and $n$
is sufficiently large (see \cite{Sen}, Thm. 1.3).
The fundamental convergence theorem (see \cite{Chu}, Thm. 1 in \S I.6 or 
\cite{Sen}, Thm. 4.2) implies that
\begin{equation}\label{fund}
\lim_{n \to \infty} q_E^{(n\df + r)}(e,f) = \df\,\la_0(f) = \frac{\df}{|E(X)|},
\end{equation}
where $\la_0$ is the unique invariant probability measure, that is,
$\la_0(f) = \frac{1}{|E|}\,.$
In view of Lemma \ref{period}, this together with (\ref{qn}) yields 
statement (b), when  $\deg(x) \ge 3$ for all $x \in X$.

Otherwise,
$$
\lim_{n \to \infty}
\frac{1}{n} \Bigl( q_E^{(1)}(e,f) + \dots + q_E^{(n)}(e,f) \Bigr) =
\frac{1}{|E(X)|}\,,
$$
and combining this with (\ref{qn}), we obtain the limit proposed in statement
(a) of Theorem \ref{limit1}.

In the case where $X$ is a cycle the $q^{(n)}(x,y)$ can be calculated
explicitly, whence the claim of the Theorem follows. This is left as
a simple exercise to the reader.

\smallskip\noindent
(c) We distinguish two cases. First, if the edge-NBRW starting at $e \in E$
is \emph{transient}, that is, the probability of returning to $e$ is
$< 1$, then $\sum_n q^{(n)}(e,f) < \infty$ for every $f \in E$, see
\cite{Chu}, Thm.\ 4 in \S I.6. Therefore, $q^{(n)}(e,f) \to 0$.

If the random walk starting at $e$ is \emph{recurrent}, i.e., it returns
to $e$ with probability $1$, then $e$ must be an essential state,
see \cite{Chu}, Thm. 4 in \S I.4 or \cite{Sen}, Lemma 5.2. 
Now by Lemma \ref{infinite},
there are infinitely many $f \in E$ such that $e \tos f$. Therefore, the
-- essential -- irreducible class $V$ of $e$ is infinite. Since the
random walk starting at $e$ does not leave $V$, we can consider the
restriction of $Q_E$ to $V$. It defines an irreducible, recurrent Markov
chain with invariant measure $\la$, the counting measure. Recurrence yields that
this is the unique invariant measure up to normalization. It has total mass
$\la(V) = \infty\,$, the chain is \emph{null recurrent}, see \cite{Chu},
\S I.6 or \cite{Sen}, \S\S 5.2--5.3. 
Therefore the convergence theorem for recurrent Markov chains yields that  
$q^{(n)}(e,f) \to 0$ for all $f \in V$. If $f \notin V$ then
$q^{(n)}(e,f) = 0$ for all $n$. 
Since $X$ is by assumption locally finite, formula (\ref{qn}) 
yields the result stated in (c).
\end{proof}

\bigskip\noindent
{\bf Uniformly irreducible random walks and amenability.} 
We now make a small detour regarding more general random walks on graphs,
recalling and improving upon the material in \cite{Wbook}, \S 10.B.

Let $X$ be a locally finite, connected graph with graph metric $d(\cdot,\cdot)$,
and  consider the transition matrix $P = \bigl(p(x,y)\bigr)_{x,y\in X}$ 
of an arbitrary 
random walk (Markov chain) on the set $X$. Then $P$ is called
\emph{uniformly irreducible} if there are constants
$K, \ep_0 > 0$ such that for any pair of neighbours $x, y$ there is
some $k \le K$ such that $p^{(k)}(x,y) \ge \ep_0$. 
Furthermore, $P$ is said to have \emph{bounded range,} if
there is $R > 0$ such that $p(x,y) > 0$ only if $d(x,y) \le R$.
These two are conditions of adaptedness of $P$ to the graph structure.

If $P$ has an \emph{invariant measure} $\nu$, then it acts on the 
Hilbert space $\ell^2(X,\nu)$ of all $F: X \to \R$ with 
$\langle F,F \rangle < \infty$, where
$
\langle F, G \rangle = \sum_x F(x)G(x)\,\nu(x)\,.
$
The operator norm satisfies $\|P\| \le 1$, and its $\ell^2$-spectral
radius is $\rho_2(P) = \lim_n \| P^n \|^{1/n}$. Note that for
$\rho(P) = \limsup_n p^{(n)}(x,y)^{1/n}$ (independent of $x,y$ by 
irreducibility) one has $\rho(P) \le \rho_2(P)$, and equality does 
not hold in general. The adjoint (more precisely, $\nu$-adjoint)
$P^*$ of $P$ on $\ell^2(X,\nu)$ has the stochastic
kernel $p^*(x,y) = \nu(y)p(y,x)/\nu(x)$.

\begin{thm}\label{improved}
Suppose that $X$ is connected, with bounded vertex degrees, and that $P$ is
uniformly irreducible with bounded range and has an invariant measure
$\nu$ satisfying $C^{-1} \le \nu(\cdot) \le C$ for some $C \ge 1$.

\smallskip
Then $\rho_2(P) = 1$ if and only if the graph $X$ is amenable.
\end{thm}

\begin{proof}[Proof (outline)] Theorem 10.6 in \cite{Wbook} states that
under the given assumptions, $\rho(P) =1$ implies amenability of $X$.
After the proof of that theorem, it is explained that the condition
$\rho(P)=1$ may be replaced with $\rho_2(P)=1$.

Conversely,  
Theorem 10.8 in \cite{Wbook} states that amenability of $X$ implies 
$\|P\|=1$. Now, let $I$ be the identity operator (or matrix), and
fix $n \ge 1$. Set $\bar P = \frac{1}{2}(I+P)$.
Then $\bar P^n$ is uniformly
irreducible, has bounded range and invariant measure $\nu$. 
If $X$ is amenable, then we get that $\|\bar P^n\| = 1$. This is
true for every $n$. Consequently, $\rho_2(\bar P) = 1$. By basic
spectral theory, also $\rho_2(P) = 1$.
\end{proof} 

More generally, the bounded range assumption can be replaced with 
tightness of the step length distributions of $P$ and $P^*$ as in 
\cite{Wbook}, Thm. 10.8.

\bigskip\noindent
We want to apply Theorem \ref{improved} not to random walks on our 
``original'' graph $X$, but to edge-NBRW on the OLG. However, the latter
is not a graph (with unoriented edges), but a digraph. Therefore, we
symmetrize it by ``removing the arrows'' from its edges. (Recall that 
the latter are ``second order'' edges, connecting edges of the original 
graph $X$). The resulting
SOLG (symmetrized oriented line graph) still has as its vertex set
the set $E$ of \emph{oriented} edges of the original graph $X$, but
neighbourhood in the SOLG is given by $e \sim f$, if $e \to f$ or $f \to e$.
We observe that in the SOLG, $q_E(e,f) > 0$ implies $e \sim f$, 
but \emph{not conversely}.

\begin{lem}\label{uniform}
If $2 \le \deg(x) \le M$ for all $x \in X$, and small cycles are dense
in $X$, then there is
$L > 0$ such that for each $e \in E$, we have $e \tos \check e$
in at most $L$ steps of edge-NBRW.\\[4pt]
In particular, $Q_E$ is uniformly 
irreducible on the symmetrized OLG.
\end{lem}

\begin{proof}
We may suppose that $X$ is infinite. Observe that the first statement
of the lemma 
implies uniform irreducibility. Indeed, let $f$ be a neighbour of $e$
in the OLG. Then either $e^+ = f^-$, in which case $q_E(e,f) \ge 1/(M-1)$,
or $f^+ = e^-$, in which case $e \tos \check e \to \check f \tos f$
in $k \le 2L+1$ steps with probability $\ge 1/(M-1)^{2L+1}$.

Now let $R > 0$ be such that $B(x,R)$ contains a cycle for every $x \in X$.
By Lemma \ref{infinite} there are infinitely many edges $f$ with $e \tos f$.
Since the vertex degree in $X$ is bounded by $M$, the number of vertices in
each $B(x,R)$ cannot exceed a certain constant $K=K(M,R)$. It follows that
$e \tos f$ for an edge $f$ with $f^+$ not contained in $B(e^-,R)$ in at most
$K$ steps of the edge-NBRW. By assumption $B(f^+,R)$ contains a cycle $C_1$ formed by edges
$e_1,\ldots,e_m$ ($m \leq K$). Since $d(e^-,f^+)>R$ neither $e$ nor
$\check e$ are edges inside the ball $B(f^+,R)$ in $X$,
and consequently neither of the two is among
the edges $e_1,\ldots,e_m$ of $C_1$.
Now, either $f \tos e_i$ (case 1) or $\check f \tos e_i$ (case 2) for some
$i\in\{1,\ldots,m\}$ in at most $R$ steps. If $f$ or $\check f=e_i$ for 
$i\in\{1,\ldots,m\}$, then 
\begin{eqnarray*} && e \tos f=e_i \to e_{i+1} \to \ldots \to e_m \to e_1 \to \ldots
    e_{i-1} \tos \check e \\ 
\mbox{or }  && e \tos f=\check e_i \to \check e_{i+1} \to \ldots \to \check e_m \to \check e_1 \to \ldots
    \check e_{i-1} \tos \check e, 
\end{eqnarray*}
respectively, in $\leq K+K+K=3K$ steps. Now let us assume that 
$f,\check f\neq e_i$ for $i\in\{1,\ldots,m\}$. Then we have in case 1 
\[ e \tos f \tos e_i \to e_{i+1} \to \ldots \to e_m \to e_1 \to \ldots
    e_{i-1} \tos \check f \tos \check e   \]
in $\leq K+R+K+R+K=2R+3K$ steps. In case 2 we have to
turn off on the way to $f$ to arrive at the cycle $C_1$. More exactly,
let $e=f_0,\ldots,f_n=f$ be a walk from $e$ to $f$ in $n\leq K$ steps.
Now consider a walk from $\check f=\check f_n$ to $e_i$ in $\leq R$ steps. It
contains some (at least one) of the edges $\check f_n, \check f_{n-1},\ldots,
\check f_{1}$. Let $\ell$ be the minimal index such that $\check f_\ell$ is not
contained in the walk. Then we have
\[ e \tos f_\ell \tos e_i \to e_{i+1} \to \ldots \to e_m \to e_1 \to \ldots e_{i-1}
      \tos \check f_\ell  \tos \check e \]
again in $\leq K+R+K+R+K=2R+3K$ steps. Thus setting $L= 2R+3K$ we have
$e \tos \check e$ in $\leq L$ steps.
\end{proof}

\begin{proof}[\bf Proof of Lemma \ref{radius}] If $X$ is not a cycle 
then by Lemma \ref{irred}, $Q_E$ is irreducible and a standard argument 
(see e.g. \cite{Sen}, Thm. 6.1, or  \cite{Wbook}, \S 1.B)
yields that 
\begin{equation}\label{rhoQ}
\rho(Q) = \rho(Q_E) = \limsup_n q_E^{(n)}(e,f)^{1/n}
\end{equation} 
is independent of
$e, f \in E$. Now apply (\ref{qn}). If $X$ is a cycle, 
$\limsup_n q^{(n)}(x,y)$
is constantly either $\frac{1}{2}$ or $1$.
\end{proof}

The fact that $\rho(Q) = \rho(Q_E)$, as stated in (\ref{rhoQ}),
is immediate from (\ref{qn}) and will be tacitly used several times.

\begin{proof}[\bf Proof of Theorem \ref{limit2}]
If $X$ is a tree then for each pair $e, f \in E$ there is at most one
$n$ such that $q_E^{(n)}(e,f) > 0$. 

Otherwise, $X$ has a cycle, and since it is almost transitive,
small cycles are dense in $X$. By Lemma \ref{radius}, $Q_E$ is
irreducible, and the OLG of $X$ is connected. Therefore the series
$\sum_n q_E^{(n)}(e,f)/\rho(Q)^n$ either converge for all $e,f$
or diverge for all $e, f \in E$, see e.g.\ \cite{Wbook}, \S 1.B.

In the convergent case, $q_E^{(n)}(e,f)/\rho(Q)^n \to 0$.

In the divergent case, edge-NBRW is \emph{$\rho$-recurrent}.
The automorphism group 
$\Ga$ of $X$ also acts with finitely many orbits on the OLG.
Therefore we can apply an adaptation of a result of {\sc Guivarc'h}
\cite{Gui}, see \cite{Wbook}, Thm.\ 7.8 and its proof: it yields
that there is a positive function $H$ on $E$ such that 
$Q_EH = \rho(Q)\cdot H$, and 
$$
q_H(e,f) = \frac{q_E(e,f)H(f)}{\rho(Q)H(e)}
$$
defines a new random walk which is $\Ga$-invariant and recurrent.
By Theorem 3.26 and Lemma 3.25 in \cite{Wbook}, $Q_H$ has an invariant
measure which is constant on each $\Ga$-orbit, and consequently
has infinite total mass. Therefore, $Q_H$ is  null recurrent,
and $q_H^{(n)}(e,f) \to 0$ for all $e, f$. Since
$$
q_H^{(n)}(e,f) = \frac{q_E^{(n)}(e,f)H(f)}{\rho(Q)^nH(e)}\,,
$$
we find that $q_E^{(n)}(e,f)/\rho(Q)^n \to 0$.
\end{proof}

A \emph{rough isometry} between two metric spaces $(X,d), (X',d')$ is a
mapping $\varphi: X \to X'$ with the following properties.
\begin{equation}\label{rough}
\begin{aligned}
A^{-1} d(x,y) - A^{-1}B &\le d'(\varphi x, \varphi y) 
                          \le A\,d(x,y) + B \quad \forall\ x, y \in X,
			  \mbox{ and}\\
d'(x',\varphi X) &\le B 
\quad \forall\ x' \in X', 
\end{aligned}
\end{equation}
where $A \ge 1$ and $B \ge 0$. 
In this case we say that the two spaces are \emph{roughly isometric}.

\begin{pro}\label{OLGri} If $X$ is a connected graph with $2 \le \deg(x) \le M$
that is not a cycle and has dense small cycles, then it is roughly isometric
with its symmetrized oriented line graph.
\end{pro}

\begin{proof} Two finite connected graphs are always roughly
isometric. Let us assume that $X$ is infinite, with edge set $E$. Throughout 
this proof, we write $d_X(\cdot,\cdot)$ for the graph distance in $X$, and 
$d_E(\cdot,\cdot)$ for the graph distance in the SOLG of $X$. Define 
the mapping $\varphi: E \to X$ by $\varphi e=e^-$.
Evidently, $\varphi$ is surjective and hence
\begin{equation}  \label{a}
d_X\bigl(x,\varphi(E)\bigr)=0 \quad \mbox{for all}\; x\in X.
\end{equation}
Now given two vertices $x,y$ in $X$ with $d_X(x,y)=d$ it is clear that
two arbitrary edges $e,f$ starting in $x$ and $y$, respectively, have
distance at least $d$ in the SOLG of $X$. It follows that
\begin{equation} \label{b}
d_X(\varphi e,\varphi f)\le d_E(e,f)\,.
\end{equation}
On the other hand, we obtain also an upper bound for $d_E(e,f)$. Clearly,
if $e,f$ are oriented the ``right way" we have $e \tos f$ in $d_X(e^-,f^-)$
steps. If one of them is oriented the other way, by Lemma \ref{uniform}
it takes at most $L$ steps to turn around, i.e.\ to reach 
$\check e$ from
$e$. Thus we have $e \tos f$ in at most $2 L + d_X(e^-,f^-)$
steps, so that
\begin{equation}  \label{c}
 d_E(e,f) - 2L \le d_X(\varphi e,\varphi f)\,.
\end{equation}
Now, setting $A=1$ and $B=2 L$ and combining (\ref{a})--(\ref{c}) yields
(\ref{rough}).
\end{proof}

\begin{proof}[\bf Proof of Proposition \ref{norm}] (a) We have 
$\| Q_E \| = \| Q_E^*Q_E \|^{1/2}$, where the adjoint operator $Q_E^*$ has 
kernel $q_E^*(e,f) = q_E(f,e)$. Let $F: E \to \R$, and let $e \in E$.
Then 
$$
Q_E^*Q_E F(e) = \sum_{f \in E} \sum_{g \in E} q_E(g,e)q_E(g,f)\,F(f)\,.
$$
Thus, $Q_E^*Q_E$ is a symmetric, stochastic operator that takes a weighted 
average of all values of $F$ on each of the finite sets 
$\{ f \in E : f^- = e^-\}$, where $e \in E$. Consequently, it has norm
$1$.\\[4pt]
(b) Instead of $Q_E$ we shall use the new transition operator 
$\bar Q_E = \frac12(I_E+Q_E)$, where $I_E$ is the identity operator.
Of course, its invariant measure is again the counting measure on $E$,
and $\bar Q_E^* = \frac12(I_E+Q_E^*)$. If we fix $n$, then 
$\bar Q_E^{*\,n} \bar Q_E^n$ is again doubly stochastic, has finite range,
and all its matrix elements are bounded below by those of $c_n Q_E$,
where $c_n = n/4^{n}$. Since $Q_E$ is (uniformly) irreducible by Lemma 
\ref{uniform}, the same holds for $\bar Q_E^{*\,n} \bar Q_E^n$.

We shall now use the obvious, but crucial relation
\begin{equation}\label{reverse}
q_E^{(n)}(e,f) = q_E^{(n)}(\check f, \check e)\,,
\end{equation}
which also holds for $\bar Q_E^n$ in the place of $Q_E^n$.
Lemma \ref{uniform} implies that for every $e \in E$,
$$
\bar q_E^{(L)}(e,\check e) \ge 1/C\,,\quad\mbox{where}\quad C=(2M)^L\,.
$$
($M$ is the upper bound on the vertex degrees.) Therefore, using
(\ref{reverse}),
$$
\bar q_E^{*\,(n)}(e,f) = \bar q_E^{(n)}(f, e)
= q_E^{(n)}(\check e, \check f) 
\le C^2\, \bar q_E^{(L)}(e, \check e)\,
\bar q_E^{(n)}(\check e, \check f)\, \bar q_E^{(L)}(\check f, f)
\le C^2 \bar q_E^{(n+2L)}(e,f)\,.
$$
In particular, we obtain that $\bar Q_E^{*\,n} \bar Q_E^n \le C^2\,
\bar Q_E^{2n+2L}$ matrix-elementwise.

Now, since $\bar Q_E^{*\,n} \bar Q_E^n$ is symmetric (self-adjoint) 
and irreducible,
Lemma 10.1 in \cite{Wbook} implies that its norm satisfies
$\| \bar Q_E^{*\,n}\bar Q_E^n\| = \rho(\bar Q_E^{*\,n} \bar Q_E^n)\,$,
the latter number being defined in the same way as in (\ref{rhoQ}), but for the
powers of $\bar Q_E^{*\,n}\bar Q_E^n$. Thus, if we take $e \in E$, then 
$$
\begin{aligned}
\rho(\bar Q_E^{*\,n} \bar Q_E^n) &= \lim_{m \to \infty}
\bigl\langle \bigl(\bar Q_E^{*\,n} \bar Q_E^n\bigr)^m\delta_e,\delta_e
\bigr\rangle^{1/m}\le 
\lim_{m \to \infty}
C^2 \, \bigl\langle \bar Q_E^{(2n+2L)m}\delta_e,\delta_e\bigr\rangle^{1/m}\\[2pt]
&= \lim_{m \to \infty} C^2 \, \bar q_E^{((2n+2L)m)}(e,e)^{1/m}
\le C^2\,\rho(\bar Q_E)^{2n+2L}\,,
\end{aligned}
$$
since $\bar q_E^{(k)}(e,e) \le \rho(\bar Q_E)^k$ for all $k \ge 0$ and 
$e \in E$, a well known fact, see e.g. \cite{Sen}, \S 6.1 or \cite{Wbook}, 
Lemma 1.9.
We infer that
$$
\rho_2(\bar Q_E) = \lim_{n \to \infty} \| \bar Q_E^{*\,n}\bar Q_E^n \|^{1/2n}
\le \lim_{n \to \infty} \bigl(C^2\,\rho(\bar Q_E)^{2n+2L}\bigr)^{1/2n} 
= \rho(\bar Q_E)\,.
$$
Since $\rho(\bar Q_E) = \frac12\bigl(1+\rho(Q_E)\bigr)$ and 
$\rho_2(\bar Q_E) = \frac12\bigl(1+\rho_2(Q_E)\bigr)\,$, we conclude that
$\rho_2(Q_E) \le \rho(Q_E)$. The reversed inequality is obvious.
\end{proof} 

\begin{proof}[\bf Proof of Theorem \ref{amenable}]
It is by now a well established fact that for connected graphs with
bounded vertex degrees, amenability is rough-isometry-invariant.
See e.g.\ \cite{Wbook}, Thm.\ 4.7 (the isoperimetric inequality $IS_{\infty}$
referred to there is the condition $\iota(X) > 0$, i.e., nonamenability),
or also the book by {\sc de la Harpe} \cite{Har}.  
Thus, in view of Proposition \ref{OLGri}, under the assumptions of
Theorem \ref{amenable} the graph $X$ is amenable if and only if its
SOLG is amenable. By (\ref{invariant}), edge-NBRW  has the counting measure 
$\la$ on $E$ as an invariant measure, and by Lemma \ref{uniform}, it is 
uniformly irreducible.
Therefore, we can apply Theorem \ref{improved} to the SOLG,
and Proposition \ref{norm}(b) allows us to replace the
$\ell^2$-spectral radius with $\rho(Q)$. 
\end{proof}

\section{Final remarks and observations}\label{final}

\noindent{\bf A.} Regarding Theorem \ref{limit1} (a+b), 
the condition  $\deg(x)\geq 3$ in Lemma \ref{period} is necessary for the
stronger convergence result of (b), as the
following example shows. Thus, if there are vertices of degree
$\leq 2$ it is in general not true that for vertex-NBRW, one has convergence of
$q^{(2n+\delta)}(x,y)$ ($\delta \in \{0,1\}$) or $q^{(n)}(x,y)$ according 
to whether $X$ is  bipartite or not (respectively).

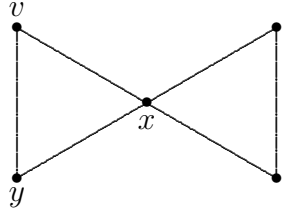
\begin{exa}
$$
\beginpicture 

\setcoordinatesystem units <.25cm,.25cm>

\setplotarea x from -5 to 5, y from -7 to -2

\plot 0 0   6.9 -4  6.9 4  0 0  -6.9 -4  -6.9 4  0 0  /
\multiput {\scriptsize $\bullet$} at 0 0  6.9 -4  6.9 4  -6.9 -4  -6.9 4 /
\put {$x$} at 0 -1
\put {$y$} at -6.9 -5
\put {$v$} at -6.9 5
\endpicture
$$
Clearly, edge-NBRW has period $\df = 3$. Write $e$ for the edge from
$y$ to $x$ and $f$ for the edge from $v$ to $y$.
We have
$$
q^{(3n)}(x,x) = 1 \AND q^{(3n+1)}(x,x) = q^{(3n+2)}(x,x) = 0 
\quad\forall\ n\,.
$$
For the edges terminating at $y$, we have $q_E^{(n)}(\check e, f) > 0$
only if $n \equiv 1 \mod 3$ and  $q_E^{(n)}(f, \check e) > 0$
only if $n \equiv 2 \mod 3$, while $q_E^{(n)}(\check e, \check e)$
and $q_E^{(n)}(f, f)$ are $>0$ only if $n \equiv 0 \mod 3$. 
Therefore, using (\ref{qn}) and (\ref{fund}), 
\[
\begin{aligned}
q^{(3n)}(y,y) &= 
\frac12\bigl( q_E^{(3n)}(\check e, \check e) + q_E^{(3n)}(f,f)\bigr)
\to \frac14\,, \\
q^{(3n+1)}(y,y) &= 
\frac12q_E^{(3n+1)}(\check e, f) \to \frac18\,,\AND\\
q^{(3n+2)}(y,y) &=
\frac12q_E^{(3n+2)}(f, \check e) \to \frac18\,, 
\quad\mbox{as $n \to \infty\,$.} 
\end{aligned} 
\]
\end{exa}

\smallskip\noindent
{\bf B.} For regular, almost transitive graphs, Lemma 3.9 of {\sc Bartholdi}
\cite{Bar}
states what is Proposition \ref{SRW} (a)+(b)+(d) and Theorems
\ref{limit1}+\ref{limit2} here. (We remark that in Lemma 3.9 of \cite{Bar}, 
the identity ``$\lim \sup_n \frac{g_n}{\beta^n}=
\lim \sup_n \frac{f_n}{\alpha^n}=\ldots$" 
should read ``$\lim \sup_n \frac{g_n}{\beta^n}=
\frac{d}{d-1}\lim \sup_n \frac{f_n}{\alpha^n}=\ldots$".)
In \cite{Bar}, a proof for SRW is suggested where one starts with the 
finite case, while for an infinite graph, one takes the sequence of balls 
$B(o,r)$ around a ``root'' vertex, applies the ``finite'' result to each ball, 
and  lets the radius tend to infinity, thereby exchanging two limits.
\cite{Bar} then suggests to use the same argument for cogrowth.
This argument has also found its way into a recent paper of
{\sc Kapovich et al.} \cite{Kap}, who state an extension to arbitrary regular graphs.
However, the argument is problematic because a priori it is by no means clear
that the two limits (for $n,r \to \infty$) may be exchanged. 

As a matter of fact, this was the starting point for the present note,
since several colleagues asked us how the mentioned argument can be made 
rigorous.
When applied to regular graphs, our method provides a simple and
rigorous proof of those statements for infinite graphs.

\smallskip\noindent
{\bf C.} Theorems \ref{limit1} and \ref{limit2} extend the corresponding
results for Cayley graphs of \cite{Wcog} to arbitrary graphs. At the same
time, the functional equation (\ref{functional}) is no more needed. The 
extension of the amenability criterion (Theorem \ref{amenable})
required more work, since the functional equation (\ref{functional})
can be used only in the regular case. Also, in the regular 
case, that criterion does not require denseness of small circles.
However, our result \emph{is} a full generalization of that amenability
criterion for (Cayley graphs of) finitely generated groups. Indeed, according
to our definition of the Cayley graph, small circles will always be dense
in the latter unless the group is freely generated by the generating set
that defines the Cayley graph. (Remember that when one of the generators
satisfies  $a_i=a_i^{-1} \ne \id$, it leads to double edges. But double 
edges give rise to circles of length $2$ according to  our definition !)

\smallskip\noindent
{\bf Acknowledgement.} The second author acknowledges discussions
with G. Noskov that stand at the origin of the questions considered in 
this paper. We also acknowledge discussions with M. Neuhauser and
a decisive hint of F. Lehner regarding the proof of Theorem \ref{amenable}.

\end{document}